\documentclass{amsart}
\usepackage{a4wide}
\usepackage[english]{babel}
\usepackage{algorithm}
\usepackage{algorithmic}
\usepackage{amsfonts}
\usepackage{amsmath}
\usepackage{amssymb}
\usepackage{amsthm}
\usepackage{bbm}
\usepackage{fourier}
\usepackage{graphicx}
\usepackage{subfigure}

\newcommand{\1}{\mathbb 1}

\newcommand{\e}{\mathrm e}

\newcommand{\diag}{\mathrm{diag}}
\newcommand{\ints}{\mathbb N}
\newcommand{\reals}{\mathbb R}

\newcommand{\KL}{\mathrm{KL}}
\renewcommand{\P}{\mathbb P}

\newcommand{\half}{\mbox{$\frac 1 2$}}
\renewcommand{\Re}{\mathrm{Re} \, }

\title{KL-learning: Online solution of Kullback-Leibler control problems}
\author{Joris Bierkens, Hilbert J. Kappen}

\begin{document}

\begin{abstract}
We introduce a stochastic approximation method for the solution of an ergodic Kullback-Leibler control problem. A Kullback-Leibler control problem is a Markov decision process on a finite state space in which the control cost is proportional to a Kullback-Leibler divergence of the controlled transition probabilities with respect to the uncontrolled transition probabilities. The algorithm discussed in this work allows for a sound theoretical analysis using the ODE method. In a numerical experiment the algorithm is shown to be comparable to the power method and the related Z-learning algorithm in terms of convergence speed. It may be used as the basis of a reinforcement learning style algorithm for Markov decision problems.
\end{abstract}

\maketitle

% \tableofcontents
% 
% \textbf{To do:} check if notation is defined at appropriate places.
% 
% \textbf{To do:} extend discussion of numerical results?

\section{Introduction}

In reinforcement learning \cite{bertsekastsitsiklis, szepesvari2010} we are interested in making optimal decisions in an uncertain environment.

Consider the setting where we are condemned to reside in a certain finite environment for an indefinite amount of time. Whenever we make a move in the environment from one state to another state, we incur a certain cost, depending on the transition. We cannot directly influence this incurred cost, but can hope to make transitions yielding a minimal average cost per transition.

This is an example of a \emph{Markov decision process} \cite{szepesvari2010} and in this paper we present a method that approximately solves this problem in a very general setting. The algorithm we present, \emph{KL-learning} (Algorithm~\ref{alg:kllearning}), observes randomly made moves (according to some Markov chain transition probabilities) and costs we incur, and finds from this, at no significant computational cost whatsoever, improved transition probabilities for the Markov chain.

This is in contrast to some other well known reinforcement learning algorithms, in which at every iteration an optimization over possible actions is necessary (e.g. Q-learning, \cite{watkins1989}) or in which an optimization step is necessary to compute optimal actions (e.g. TD-learning, \cite{sutton1988}).

The background for this method is the setting of the \emph{Kullback-Leibler (KL) control problem}, introduced in \cite{todorov2007}. A KL-control problem is a Markov decision process in which the control costs are proportional to a \emph{Kullback-Leibler divergence} or \emph{relative entropy}. In \cite{todorov2007} also a reinforcement style learning algorithm (Z-learning) was presented, which operates under the assumption of there being an absorbing state in which no further costs are incurred. This assumption is not made in our algorithm; instead we assume ergodicity of the underlying Markov chain. Arguably this yields a more general setting, in which a hard reset of the algorithm is never necesssary. KL control problems may also be solved using techniques from graphical model inference \cite{KappenOpper:2012}.

As a preliminary, we introduce the KL-control setting in Section~\ref{sec:klcontrol}.
In Section~\ref{sec:kllearning} the KL-learning algorithm is presented and motivated on a heuristic level. We then describe the ODE method \cite{benvenisteetal, kushnerclark, kushneryin, ljung} in Section~\ref{sec:odemethod} along with an application to a stochastic gradient algorithm and Z-learning \cite{todorov2007} as illustrative examples. We then apply the ODE method to KL-learning in Section~\ref{sec:convergence}. A numerical example is provided in Section~\ref{sec:numerical} after which a short discussion follows in Section~\ref{sec:discussion}.

\section{Kullback-Leibler control problems}
\label{sec:klcontrol}
In this section we introduce the particular form of Markov decision process which have a particularly convenient solution. We will refer to these problems as Kullback-Leibler control problems. For a more detailed introduction, see \cite{todorov2007}.

Let $t = 0, 1, 2, \hdots$ denote time. Consider a Markov chain $(X_t)_{t=0}^{\infty}$ on a finite state space $S = \{1, \hdots, n \}$ with transition probabilities $q = [q(j|i)]$ which we call the \emph{uncontrolled dynamics}. We will make no distinction between the notation $q_{ij}$ and $q(j|i)$, where $q(j|i) = q_{ij}$ denotes the probability of jumping from state $i$ to state $j$; the notation $q_{ij}$ will me more convenient when working with matrices. 

Suppose for every jump of the Markov chain from state $i$ to state $j$ in $S$ a \emph{transition dependent cost} $c(j|i)$ is incurred. Sometimes we will use the notation $c(i)$ to denote costs depending only state, i.e. $c(j|i) = c(i)$ for all $i, j = 1, \hdots, n$. A state $i$ is called \emph{absorbing} if $q(i|i) = 1$ and $c(i)=0$. 

We wish to change the transition probabilities in such a way as to minimize for example the total incurred cost (assuming there exist absorbing states where no further costs are incurred) or the average cost per stage. For deviating from the transition probabilities control costs are incurred equal to
\[ \frac 1 {\beta} \KL(p(X_{t+1} | X_t) || q(X_{t+1} | X_t)) = \frac 1 {\beta} \sum_{j=1}^n p(X_{t+1} = j | X_t) \ln \left( \frac{p(X_{t+1} = j | X_t)}{q(X_{t+1} = j | X_t)} \right)\] at every time step, in addition to the cost per transition $c(X_t|X_{t-1})$, where $\beta > 0$ is a weighing factor, indicating the relative importance of the control costs.

To put this problem in the usual form of a discrete time stochastic optimal control problem, we write $p_{ij} = \exp(u_j(i)) q_{ij}$. This guarantees positive probabilities and absolute continuity of the controlled dynamics with respect to the uncontrolled dynamics. In the case of an infinite horizon problem and minimization of a total expected cost problem, the corresponding Bellman equation for the value function $\Phi$ is
\[ \Phi(i) = \min_{(u_1, \hdots, u_n) \in \reals^n} \left\{\sum_{j=1}^n c(j|i) + q(j|i) \exp(u_j) ( u_j / \beta + \Phi(j) ) \right\},\]
where the minimization is over all $u_1, \hdots u_n$ such that $\sum_{j=1}^n \exp(u_j) q(j|i) = 1$.
If there are no absorbing states, the total cost will always be infinite and the expression above has no meaning. We may then instead aim to minimize the expected average cost.
For an average cost problem, the Bellman equation for the value function $\Phi$ is
\begin{equation} \label{eq:bellmanaveragecost} \rho + \Phi(i) = \min_{(u_1, \hdots, u_n) \in \reals^n} \left\{ \sum_{j=1}^n c(j|i) +  q(j|i) \exp(u_j) (u_j / \beta + \Phi(j)) \right\},\end{equation}
where again the minimization is over all $u_1, \hdots u_n$ such that $\sum_{j=1}^n \exp(u_j) q_{ij} = 1$, and where $\rho$ is the optimal average cost.
In the average cost case we restrict the possible solutions by requiring that
\begin{equation}
\label{eq:restrictionphiaveragecost}
\sum_{i=1}^n \exp(-\beta \Phi(i)) = 1;
\end{equation}
otherwise any addition by a scalar would result in another possible value function. The reason for the particular form of this restriction will become clear later. 

Note that in case the total expected cost problem has a finite value function, the solution of the average cost problem~\eqref{eq:bellmanaveragecost} would have a solution with $\rho = 0$. This shows that in a sense the average cost problem is more general, since then~\eqref{eq:bellmanaveragecost} remains valid for the total expected cost problem. Therefore we will henceforth only consider the average cost problem case. 

So far the derivations have been standard; see \cite{bertsekas1} for more information on dynamic programming and the Bellman equation.

It is remarkable that a straightforward computation using Lagrange multipliers, as in \cite{todorov2007}, yields that the optimal $u_j(i)$ and value function $\Phi$ solving~\eqref{eq:bellmanaveragecost} are given by the simple expressions
\begin{equation} \label{eq:optimalsolutionaveragecost} u_j^*(i) = \ln(z_j^* / \lambda^* z_i^*) - \beta c(j|i)), \quad  \quad \Phi(i) = - \frac 1 {\beta} \ln ( z_i^*), \end{equation}
with $z^* \in \reals^n$ given implicitly by
\[ \lambda^* z_i^* = \sum_{j=1}^n \exp(-\beta c(j|i))q(j|i) z_j^*,\]
which may be written as $\lambda^* z^* = H z^*$, with 
\begin{equation} \label{eq:definitionH} H = [h_{ij}] \quad \mbox{with entries} \quad h_{ij} = \exp(-\beta c(j|i) )q(j|i)
\end{equation}
and where $\lambda^* = \exp(- \beta \rho^*)$. This $z^*$ should be normalized in such a way that the value function agrees with the value 0 in the absorbing states for a total expected cost problem, or with the normalization~\eqref{eq:restrictionphiaveragecost} in the average cost case, which is chosen in such a way that it corresponds to $||z^*||_1 = \sum_{i=1}^n z_i^* = 1$. 
The optimal transition probabilities simplify to
\[ p(j|i)^* = q(j|i) \exp(-\beta c(j|i)) \frac{z_j^*}{\lambda^* z_i^*}.\]

According to Perron-Frobenius theory of non-negative matrices (see \cite{hornjohnsonmatrixanalysis}), if the uncontrolled Markov chain $q$ is irreducible then there exists, by Observation~\ref{obs:Hirreducible} below, a simple eigenvalue $\lambda^*$ of $H$ equal to the spectral radius $\rho(H)$, with an eigenvector $z^*$ which has only positive entries. Since $\lambda^*$ is a simple eigenvalue, $z^*$ is unique op to multiplication by a positive scalar. These $\lambda^*$ and $z^*$ (with $z^*$ normalized as above) are called the \emph{Perron-Frobenius} eigenvalue and eigenvector, respectively. The optimal average cost is given by $\rho^* = - \frac 1 {\beta} \ln \lambda^*$. In case of a total expected cost problem, where $\rho^* = 0$, it follows that $\lambda^* = 1$, which may also be shown directly by analysis of the matrix $H$. 

Recall that a nonnegative matrix $A$ is called \emph{irreducible} if for every pair $i, j \in S$, there exists an $m \in \ints$ such that $(A^m)_{ij} > 0$. In particular, a Markov chain $p$ is called irreducible if the above property holds for its transition matrix.

\subsection{Observation}
\label{obs:Hirreducible}
Suppose the finite Markov chain $q$ on $S = \{ 1, \hdots, n\}$ with transition probabilities $q(j|i)$ is irreducible. Then $H$ as given by~\eqref{eq:definitionH} is irreducible.
In particular, there exists a unique (modulo scalar multiples) positive eigenvector $z^* \in \reals^n$ of $H$ such that $H z^* = \lambda^* z^*$, where $\lambda^* = \rho(H) = \sup_{\mu \in \sigma(H)} |\mu|$, the \emph{spectral radius} of $H$.

\emph{Proof.}
Let $\gamma = \min_{i \in S} \e^{-\beta c(j|i)}$. Let $i,j \in S$ and pick $m \in \ints$ such that $(q^m)_{ij} > 0$.
Then 
\[ (H^m)_{ij} = \sum_{k_1=1}^n \hdots \sum_{k_{m-1} = 1}^n H_{i k_1} H_{k_1 k_2} \hdots H_{k_{m-1}, j} 
\geq \gamma^m \sum_{k_1=1}^n \hdots \sum_{k_{m-1} = 1}^n q_{i k_1} q_{k_1 k_2} \hdots q_{k_{m-1}, j} > 0.\]

The existence and uniqueness of the eigenvalue and corresponding eigenvector is then an immediate corollary of the Perron-Frobenius theorem \cite[Theorem 8.4.4]{hornjohnsonmatrixanalysis}.
\qed

Recall that a Markov chain $[p_{ij}]$ is said to satisfy detailed balance if there exists a probability distribution $(p_i)$ such that $p_i p_{ij} = p_j p_{ji}$ for all $i, j$. In this case $(p_i)$ is an invariant probability distribution for the Markov chain.

\subsection{Proposition}
Suppose the uncontrolled dynamics $q$ satisfy detailed balance (with respect to the invariant probability distribution given by $(q_i)$).

\begin{itemize}
\item [(a)] If the transition costs are actually state costs, i.e. $c(j|i) = c(i)$ for $i,j=1, \hdots, n$, then the optimal controlled dynamics satisfy detailed balance with invariant probability distribution given by 
\[ p_i \propto q_i \exp(\beta c(i)) (z_i^*)^2, \quad i =1, \hdots, n.\]
\item[(b)] If the transition costs are symmetric, i.e. $c(j|i) = c(i|j)$ for $i,j =1,\hdots, n$, then the optimal controlled dynamics satisfy detailed balance with invariant probability distribution give by
\[ p_i \propto q_i (z_i^*)^2, \quad i =1, \hdots, n.\]
\end{itemize}

\emph{Proof.}
We will prove (a), the proof of (b) is analogous.

Using that $p_{ij} = \exp(u_j^*(i)) q_{ij}$ with $u_j^*(i)$ given by~\eqref{eq:optimalsolutionaveragecost}, we verify that $p_i p_{ij} = p_j p_{ji}$ for all $i, j$. Indeed, 
\begin{align*}
p_i p_{ij} & = q_i \exp(\beta c(i)) (z_i^*)^2 q_{ij} z_j^*/z_i^* \exp(- \beta c(i)) / Z= q_i q_{ij} z_i^* z_j^* / Z \\
& = q_j q_{ji} z_i^* z_j^* / Z = q_j \exp(\beta c(j)) (z_j^*)^2 q_{ji} z_i^*/z_j^* \exp(- \beta c(j) / Z = p_j p_{ji},
\end{align*}
where $Z = \sum_{k=1}^n q_k \exp(\beta c(k)) (z_k^*)^2$ is a normalization constant.\qed

\subsection{Example: solution in case of trivial detailed balance}
If we take as uncontrolled dynamics $q_{ij} = q_j$, where $q_j$ is a probability distribution on $\{1, \hdots, n\}$, then $H_{ij} = \exp(- \beta c(i)) q_j$ is of rank one and has non-zero eigenvalue $\lambda^* = \sum_{j=1}^n q_j \exp(-\beta c(j))$ with eigenvector $z^*$ given by $z^*_i =\exp(-\beta c(i))$. The optimal transition probabilities are given by
\[ p_{ij} = q_j \exp(-\beta c(j)) / \sum_{k=1}^n q_k \exp(- \beta c(k)),\]
which again are independent of $i$. Therefore the Markov chain given by the controlled dynamics has invariant probability distribution $[p_j] = [p_{ij}]$. The optimal average cost is given by
\[ \rho^* = -\frac 1{\beta} \ln \lambda^* = - \frac 1{\beta}\ln \sum_{k=1}^n q_k \exp(-\beta c(k)).\]
\hfill $\diamond$

\section{KL-learning}
\label{sec:kllearning}
As explained in the previous section, a Kullback-Leibler control problem may be solved by finding the Perron-Frobenius eigenvalue $\lambda^*$ and eigenvector $z^*$ of the matrix $H$ given by~\eqref{eq:definitionH}.

A straightforward way to find $\lambda^*$ and $z^*$ is using the \emph{power method}, i.e. by performing the iteration
\begin{equation} \label{eq:powermethod} z_{k+1} = \frac{ H z_k}{||H z_k||}.\end{equation}
This assumes that we have access to the full matrix $H$. Our goal is to relax this assumption, and to find $z$ by iteratively stepping through states of the Markov chain using the uncontrolled dynamics $q$, using only the observations of the cost $c(j|i)$ when we make a transition from state $i$ to state $j$..

In \cite{todorov2007} a \emph{stochastic approximation algorithm} (see \cite{benvenisteetal, bertsekastsitsiklis, kushnerclark, kushneryin}), referred to as Z-learning, is introduced for the case $\lambda^* = 1$. We will extend this method here to the case where $\lambda^*$ is a priori unknown.

In this section we will denote vectors by bold letters, e.g. $v$. Components of this vector will be denoted as $v(i)$ or $v_i$. The notation $\1$ is used for the column vector containing only ones. A vector $v \in \reals^n$ is said to be nonnegative ($v \geq 0$) if $v(i) \geq 0$ for all $i = 1,\hdots, n$ and positive ($v > 0$) if $v(i) > 0$ for all $i = 1, \hdots, n$.

The algorithm we will consider is Algorithm~\ref{alg:kllearning}. The parameter $M \in \ints$ denotes the number of iterations of the algorithm, and $\gamma > 0$ indicates the stepsize. We assume that the Markov transition probabilities $q(\cdot|\cdot)$ are irreducible and aperiodic, and hence ergodic.

\begin{algorithm}[ht]
\caption{KL-learning}
\label{alg:kllearning}

\begin{algorithmic}
\STATE $z \gets \frac 1 n \1$, $\lambda \gets 1$, $x \gets \mbox{any state in $S$}$
\FOR{$k = 1 \ \TO \ M$}
\STATE $y \gets \mbox{independent draw from $q(\cdot|x)$}$
\STATE $\Delta \gets \exp(-\beta c(y|x)) z(y)/ \lambda - z(x)$
\STATE $z(x) \gets z(x) + \gamma \Delta $
\STATE $\lambda \gets \lambda + \gamma \Delta$
\STATE $x \gets y$
\ENDFOR
\end{algorithmic}

\end{algorithm}

At every iteration, we make a random jump to a new state. Based on our observation of the incurred cost at the previous step, and current values of $\lambda$ and two components of $z$, a number $\Delta$ is computed that says how much $z$ and $\lambda$ should be changed. The value of $\lambda$ is always equal to $\sum_{i=1}^n z_i = ||z_i||_1$. Note that every step of the iteration consists of only simple algebraic operations and hence has time complexity $\mathcal O(1)$. In particular, no optimization is needed, as opposed to e.g. $Q$-learning \cite{watkins1989}.

A theoretical analysis of (a slightly modified version of) this algorithm will be performed in Section~\ref{sec:convergence}. The results of that section are summarized in Theorem~\ref{thm:convergencekllearning}. First we provide some intuition.

\subsection{Heuristic motivation}
Suppose at time $m$ we are in state $i$. The expected value of $\Delta$ is 
\[ \sum_{j=1}^n q_{ij} \left( \exp(-\beta c_{ij}) z(j)/ \lambda - z(i) \right) = (H z)_i / \lambda - z_i.\]
Since $\lambda = ||z||_1$, the update to $z$ may be interpreted as 
\[ z_{\mathrm{new}} = z + \gamma \left[ (H z)(i) / \lambda - z(i) \right] = (1 - \gamma) z(i) + \gamma (Hz)(i) / ||z||_1,\]
a convex combination of the old value of $z(i)$ and the value $z(i)$ would obtain after an iteration of the power method described above. The normalization is however based on the previous value of $z$ but this does not affect the convergence of the power method.

The frequency of updates to the $i$-th component of $z$ depends, on the long run, on the equilibrium distribution $(q_i)$ of the underlying Markov chain. This will be a major concern in the convergence analysis of the algorithm. It will turn out that the convergence of the algorithm will depend on the stability properties of a certain matrix, $A$ say. If we wish the algorithm to converge for a certain invariant distribution, this corresponds to the matrix $D A$ being stable, where $D$ is a diagonal matrix with the invariant distribution on the diagonal. This will be made clear in Section~\ref{sec:convergence}.

\section{Analysis of stochastic approximation algorithms through the ODE method}
\label{sec:odemethod}

In this section a general and powerful method for analyzing the behaviour and possible convergence of stochastic approximation algorithms is described. It will be applied to Algorithm~\ref{alg:kllearning} in Section~\ref{sec:convergence}. This method, called the \emph{ODE method}\footnote{Here ODE is an abbreviation for \emph{ordinary differential equation}.}, was first introduced by Ljung \cite{ljung} and developed significantly by Kushner and coworkers \cite{kushnerclark, kushneryin}. Accounts that are well suited for computer scientists and engineers may be found in \cite{benvenisteetal, bertsekastsitsiklis}.

The theory is illustrated by applying it to some stochastic algorithms. The new contribution of this section to the existing theory is the necessity of \emph{diagonal stability} for the convergence of certain stochastic algorithms, as discussed in Section~\ref{sec:diagonalstabilitynecessary}.

\subsection{Outline of the ODE method}
The idea of the ODE method is to establish a relation between the trajectories of a stochastic algorithm with decreasing stepsize, and the trajectories of an ordinary differential equation. If all trajectories of the ODE converge to a certain equilibrium point, the same can then be said about trajectories of the stochastic algorithm. This is made more precise in the following theorem, which is a special case of \cite[Theorem 6.6.1]{kushneryin} tailored to our needs.

\subsection{Hypotheses}
\label{hyp:generalconvergence}
Consider the general stochastic approximation algorithm given by Algorithm~\ref{alg:general}, assuming the following assumptions and notation:

\begin{algorithm}[ht]
\caption{General stochastic approximation algorithm for theoretical analysis}
\label{alg:general}
\begin{algorithmic}
\STATE $x_0 \gets \mbox{any state in} \ S$
\STATE $\theta_0 \gets \mbox{any initial vector in $\reals^n$}$
\FOR{$k = 1 \ \TO \ \infty$}
\STATE $x_k \gets \mbox{independent fraw from $q(\cdot|x_{k-1})$}$
\STATE $\theta_k \gets \theta_{k-1} + \gamma_k f(\theta_{k-1}, x_{k-1} ,x_k)$
\ENDFOR
\end{algorithmic}
\end{algorithm}

\begin{itemize}
\item[(i)]
Let $\gamma_1, \gamma_2, \hdots$ be a sequence of step sizes, satisfying $\sum_{k=1}^{\infty} \gamma_k = \infty$ and $\sum_{k=1}^{\infty} \gamma_k^2 < \infty$;
\item[(ii)]
Let $q(\cdot|\cdot)$ be irreducible aperiodic Markov transition probabilities on a finite state space $S$ with invariant probabilities $q_i$, $i \in S$; 
\item[(iii)]
Suppose that $\{ \theta_k : k \in \ints\} \subset K$ with probability one, where $K$ is some compact (i.e. closed and bounded) subset of $\reals^n$;
\item[(iv)]
Suppose $(\theta, x, y) \mapsto f(\theta, x, y) : K \times S \times S \rightarrow \reals^n$ is continuous in $\theta$ for every $x, y \in S$.
\item[(v)] Define $\overline g : K \rightarrow \reals^n$ by
\begin{equation}
\label{eq:averageg} \overline g(\theta) := \sum_{i \in S} \sum_{j \in S} q_i q_{ij} f(\theta, i, j).
\end{equation}
\item[(vi)] Define $t_0 := 0$ and $t_k := \sum_{i=1}^k \gamma_i$ for $k \in \ints$.
% For $t \geq 0$ let $m(t)$ denote the unique value of $k$ such that $t_k \leq t < t_{k+1}$, and for $t < 0$ set $m(t) = 0$. 
Denote, for all $t \geq 0$ and $k = 0, 1, 2, \hdots$, $\theta^k(t) := \theta_p$ 
% and $Z^k(t) := Z_p$ 
for the unique $p$ such that $t_p \leq t + t_k < t_{p + 1}$, and $\theta^k(t) := 0$ if $t + t_k < t_0$.

% \item[(iii)]
% Let $H$ denote a bounded and closed subset of $\reals^n$ of the form $H = \{ \theta \in \reals^n : a_i \leq \theta^i \leq b_i \}$, with $-\infty < a_i < b_i < \infty$ for $i = 1, \hdots, n$;
% \item[(iv)] Let $\Pi_H: \reals^n \rightarrow H$, $\Pi_H(x) := \arg \min_{y \in H} |x - y|$ denote the projection on $H$ (where $|\cdot|$ denotes the Euclidean norm on $\reals^n$);
% \item[(vi)] Let $Z_k$ be random variables satisfying $\theta_k = \theta_{k-1} + \gamma_k g(\theta_{k-1}, x_{k-1}, x_k) + \gamma_k Z_k$ for $k \in \ints$, i.e. $Z_k$ is the correction to project the $k$-th update to $\theta$ back into $H$;

\end{itemize}

These assumptions are sufficient for our purposes. The sequence $\gamma$ denotes the \emph{stepsize} or \emph{gain}. The conditions under (i) on $\gamma$ are standard conditions to guarantee that the gain gradually decreases, but not too quickly, in which case the algorithm would stop making significant updates before being able to converge.

In \cite{kushneryin} more general classes of algorithms and assumptions are considered.

\subsection{Theorem (convergence of stochastic algorithms with state dependent updates)}
\label{thm:convergencestochasticalgorithm}
Suppose Assumptions~\ref{hyp:generalconvergence} hold. Then, with full probability,
\begin{itemize}
\item[(i)] Every sequence in the collection of functions $\{ \theta^k : k \in \ints \}$ (as defined under Assumption~\ref{hyp:generalconvergence} (vi)) admits a convergent subsequence with a continuous limit;\footnote{Here by convergence we mean uniform convergence on bounded intervals.}
\item[(ii)] Let $\theta$ denote the limit of some converging subsequence in $\{ \theta^k : k \in \ints \}$ (which always exists by (i)). Then $\theta$ satisfies the ODE
\begin{equation}
\label{eq:ODE}
\dot \theta = \overline g(\theta)
%  + z,
\end{equation}
% where $z = z(\theta) \in \reals^n$ is the minimum force necessary to keep the trajectory of~\eqref{eq:ODE} in $H$ (see \cite[Section 4.3]{kushneryin}).
\item[(iii)] If a set $A \subset \reals^n$ is globally asymptotically stable with respect to the ODE~\eqref{eq:ODE}, then $\theta_k \rightarrow A$, i.e. $\min_{x \in A} |\theta_k - x| \rightarrow 0$.
\end{itemize}

\emph{Outline of proof.}
The proof consists of a verification of the conditions of \cite[Theorem 6.6.1]{kushneryin}. One key ingredient for this verification is Lemma~\ref{lem:convergencetoinvariantdistribution} below, which says that convergence of the pair $(x_{k-1}, x_k)$ to its equilibrium distribution $(q_i q_{ij})_{i, j \in S}$ happens exponentially fast.\hfill $\diamond$

% \textbf{To do/consider:} more extensive outline of proof? / intuition of the theorem

% Let $\widetilde \E$ denote expectation under this invariant distribution for $(x_{k-1}, x_k)$. Then 
% \[ \frac 1 { \gamma_k } \widetilde \E [ \theta_k - \theta_{k-1} ] = \widetilde \E f(\theta_{k-1}, x_{k-1}, x_k) = \sum_{i \in S} \sum_{j \in S} q_i q_{ij} f(\theta_{k-1}, i , j) = \overline g(\theta_{k-1}).
% \]
% This explains the specific form ~\eqref{eq:ODE} of the differential equation.\qed
% 
Recall the \emph{total variation distance} \cite[Section 4.1]{levinetal} of two probability measures $\mu_1, \mu_2$ on a discrete space $S$,
\[ || \mu_1 - \mu_2 ||_{\mathrm{TV}} := \sup_{A \subset S} |\mu_1 (A) - \mu_2(A)| = \half \sum_{i \in S} | \mu_1(i) - \mu_2(i) |.\]

\subsection{Lemma (Markov chain convergence to invariant distribution)}
\label{lem:convergencetoinvariantdistribution}
Let $q(j|i)$, $i, j \in S$, denote the transition probabilities of an irreducible, aperiodic Markov chain $X$ on a finite state space $S$ with invariant distribution $q_i$, $i \in S$.
Let $\mu_k^x$ be the probability measure on $S \times S$ denoting the distribution of $(X_{k-1}, X_k)$ given $X_0 = x$. Let $\overline{\mu}$ denote the probability measure on $S \times S$ given by $\overline{\mu}(i, j) = q_i q(j|i)$.

Then there exist constants $\alpha \in (0, 1)$ and $C > 0$ such that
\begin{equation}\label{eq:convergencetoinvariantdistribution} \max_{x \in S} || \mu_k^x - \overline{\mu} ||_{\mathrm{TV}} \leq C \alpha^k, \quad \mbox{for all} \quad k \in \ints.\end{equation}

\emph{Proof:}
Let $\nu_k^x$ denote the probability measure on $S$ denoting the distribution of $X_k$ given initial condition $X_0 = x$. By~\cite[Theorem 4.9]{levinetal}, there exist constants $\widetilde C > 0$ and $\alpha \in (0,1)$ such that
\[ \max_{x \in S} ||\nu_{k-1}^x - q ||_{\mathrm{TV}} \leq \widetilde C \alpha^{k-1} \quad \mbox{for} \ k \in \ints.\]
Therefore
\begin{align*}
 \max_{x \in S} || \mu_k^x - \overline{\mu} ||_{\mathrm{TV}} & = \max_{x \in S} \half \sum_{i \in S} \sum_{j \in S} | \P(X_{k-1} = i, X_k = j | X_0 = x) - q_i q(j|i)| \\
& = \max_{x \in S} \half \sum_{i \in S} q(j|i) \sum_{j \in S} | \P(X_{k-1} = i | X_0 = x) - q_i  |  \\
& = \max_{x \in S} \half \sum_{i \in S} | \P(X_{k-1} = i | X_0 = x) - q_i | = \max_{x \in S} ||\nu_{k-1}^x - q ||_{\mathrm{TV}} \leq \widetilde C \alpha^{k-1}.
\end{align*}
By letting $C = \widetilde C / \alpha$ we find that ~\eqref{eq:convergencetoinvariantdistribution} holds.\qed

\subsection{Remark}
Note that a boundedness assumption is made in Theorem~\ref{thm:convergencestochasticalgorithm}. In practice, this is not an unreasonable assumption, since float sizes are bounded in many programming languages. The boundedness may be enforced by a projection step in the algorithm, leading to a slightly more complex formulation of Theorem~\ref{thm:convergencestochasticalgorithm}. In particular, the resulting ODE becomes a projected ODE. See \cite[Section 4.3]{kushneryin}.

\subsection{Example: A stochastic gradient algorithm}
Suppose we wish to minimize a function $h : \reals^n \rightarrow \reals$ with bounded first derivatives but that we do not have full access to its gradient $g_j = \frac{\partial h}{\partial \theta^j}$. Instead, the observations we make are determined by an underlying Markov chain $(x_k)$ on the state space $\{1, \hdots,n\}$ with aperiodic, irreducible transition probabilities $q_{ij}$. In case a jump is made to $x_k$, we observe $\frac{\partial h}{\partial \theta^{(x_k)}}(\theta)$ for some $\theta \in \reals^n$. Is there a stochastic approximation algorithm that can minimize $h$ under these restrictive conditions?

Consider Algorithm~\ref{alg:stochasticgradient}. We use $e_i$ to denote the unit vector in direction $i$.

\begin{algorithm}
\caption{Stochastic gradient algorithm}
\label{alg:stochasticgradient}
\begin{algorithmic}
\STATE $x_0 \gets \mbox{any state in} \ S$
\STATE $\theta_0 \gets \mbox{any initial vector in $\reals^n$}$
\FOR{$k = 1 \ \TO \ \infty$}
\STATE $x_k \gets \mbox{independent fraw from $q(\cdot|x_{k-1})$}$
\STATE $\theta_k \gets \theta_{k-1} - \gamma_k \frac{\partial h(\theta_{k-1})}{\partial \theta^{(x_k)}} e_{(x_{k-1})}$
\ENDFOR
\end{algorithmic}

\end{algorithm}

Since $\nabla h$ is bounded, the trajectories of this algorithm are restricted to the bounded set 
\[ K = \left\{ \theta \in \reals^n : |\theta | \leq \max(|\theta_0|, |\nabla h|) \right\}\] 
with probability one.
The corresponding ODE (in the sense of Theorem~\ref{thm:convergencestochasticalgorithm}) is~\eqref{eq:ODE} with
\begin{equation} \overline g^i(\theta) = - q_i \sum_{j=1}^n q_{ij} \frac{\partial h(\theta)}{\partial \theta^j}.\end{equation}
Let $R = [r_{ij}]$ be the matrix defined by $r_{ij} = q_i q_{ij}$, $i, j = 1, \hdots, n$. We may then write
\begin{equation} \label{eq:stochgradientODE} \dot g(\theta) = - R \nabla h(\theta).\end{equation}

Clearly the minimum $\theta^*$ of $h$, where $\nabla h(\theta^*) = 0$, gives an equilibrium point of this ODE. It is not immediately clear whether this is the only equilibrium point.

We will now make the following assumptions: 
\begin{align}
\label{hyp:Rposdef}
\mbox{The matrix $R$ given by the entries $r_{ij} := q_i q_{ij}$ is symmetric, positive definite.}  \\
\label{hyp:hconvex}
\mbox{The function $h$ is twice differentiable and strictly convex.}
\end{align}
Sufficient conditions for~\eqref{hyp:Rposdef} to hold are the following:
\begin{itemize}
\item[(i)] The Markov chain given by $q_{ij}$ satisfies \emph{detailed balance}, i.e. $q_i q_{ij} = q_j q_{ji}$ for $i, j = 1, \hdots, n$;
\item[(ii)] The Markov chain given by $q_{ij}$ is \emph{strictly lazy}, i.e. $q_{ii} > \half$ for $i =1, \hdots, n$.
\end{itemize}
Indeed, if these conditions are satisfied, then $R$ is symmetric by the detailed balance condition. Since the Markov chain is lazy, $R$ is strictly row diagonally dominant, so that all its eigenvalues are positive.

Under Assumption~\eqref{hyp:Rposdef}, define a Lyapunov function $V : \reals^n \rightarrow \reals$ by
\[ V(\theta) = \half \langle \nabla h(\theta), R \nabla h(\theta) \rangle,\]
where $\langle \cdot, \cdot \rangle$ denotes the Euclidean inner product on $\reals^n$.

Write $H(\theta) = \left(\frac{\partial^2 h}{\partial \theta^i \theta^j}\right)_{i,j=1,\hdots, n}$ to denote the Hessian matrix of $h$ at $\theta$, and note that, since $h$ is strictly convex, the matrix $H(\theta)$ is positive definite for all $\theta \in \reals^n$.
Then if $\theta(t)$ satisfies~\eqref{eq:stochgradientODE},
\begin{align*}
\frac{d}{dt} V(\theta(t)) & = \langle \nabla h(\theta(t)), R H(\theta(t)) \dot \theta(t) \rangle = - \langle \nabla h(\theta(t)), RH(\theta(t)) R \nabla h(\theta(t)) \rangle \leq 0,
\end{align*}
with strict inequality if $\nabla h(\theta) \neq 0$. This shows that the ODE~\eqref{eq:stochgradientODE} is globally asymptotically stable with unique equilibrium $\theta^*$ satisfying $\nabla h(\theta^*) = 0$. By Theorem~\ref{thm:convergencestochasticalgorithm} (iii) therefore Algorithm~\ref{alg:stochasticgradient} converges almost surely to $\theta^*$.

In this case we are in some sense lucky to be able to find a Lyapunov function to establish global stability of the ODE. In the case of Algorithm~\ref{alg:kllearning} (KL-learning) we have not yet found a global Lyapunov function and so far can only achieve local stability around the equilibrium in certain cases. For illustrative purposes, we now also perform such a local analysis to the current example.

It is immediately clear that under assumption~\eqref{hyp:Rposdef}, the only equilibrium of the ODE~\eqref{eq:stochgradientODE} satisfies $\nabla h(\theta) = 0$. It remains to establish the stability of the ODE around that equilibrium point. The linearized version of the ODE around the equilibrium is
\begin{equation}
\label{eq:linearizedODEstochasticgradient}
\dot \theta(t) = - R H(\theta^*) \theta.
\end{equation}
We therefore need to determine the spectrum of the matrix $R H(\theta^*)$. Indeed $R H(\theta^*) R + R H(\theta^*) R$ is positive definite, so that by Lyapunov's theorem \cite[Theorem 2.2.1]{hornjohnsontopics} $R H(\theta^*)$ has only eigenvalues in the open right halfplane. We may conclude from this local analysis, by the Hartman-Grobman theorem \cite[Section 2.8]{perko}, that the equilibrium $\theta^*$ is locally asymptotically stable.

\subsection{Remark.}
Under the assumption that we can only observe $\frac{\partial h}{\partial \theta^j}$ if we jump to state $j$, a simpler algorithm would consist of the update rule $\theta_k \gets \theta_{k-1} - \gamma_k \frac{\partial h (\theta_{k-1})}{\partial \theta^{(x_k)}} e_k$, i.e. to update the $(x_k)$-th component of $\theta$ instead of the $(x_{k-1})$-th component. In this case a Lyapunov function would be given by $V(\theta) = \sum_{i=1}^n q_i \left(\frac{\partial h(\theta)}{\partial \theta^i}\right)^2$ and Assumption~\eqref{hyp:Rposdef} would not be required. However, the analysis of Algorithm~\ref{alg:stochasticgradient} has more in common with the upcoming analysis of Algorithm~\ref{alg:kllearning} (KL-learning), because in that algorithm the updates also depend upon the previous and current state of the Markov chain.

\subsection{Example: Z-learning}
\label{sec:zlearning}

In \cite{todorov2007}, the Z-learning algorithm is presented as a way to solve the eigenvector problem $H z^* = z^*$, where $H = [h_{ij}]$ is a nonnegative irreducible matrix with spectral radius $\rho(H) = 1$ of the form $h_{ij} = \exp(-\beta c_{ij}) q_{ij}$ as in Section~\ref{sec:klcontrol}, with $[q_{ij}]$ the transition probabilities of some irreducible Markov chain on $S = \{ 1, \hdots, n\}$. This problem is an important special case of the problem we address in this paper, namely solving $H z^* = \lambda^* z^*$ with unknown spectral radius $\rho(H) = \lambda^*$.

The Z-learning algorithm is given by Algorithm~\ref{alg:zlearning}.

\begin{algorithm}[ht]
\caption{Z-learning}
\label{alg:zlearning}

\begin{algorithmic}
\STATE $z_0 \gets \1$, $x_0 \gets \mbox{any state in $S$}$
\FOR{$k = 1 \ \TO \ \infty$}
\STATE $x_k \gets \mbox{independent draw from $q(\cdot|x_{k-1})$}$
\STATE $z_k \gets z_{k-1} + \gamma_k \left( \exp(-\beta c(x_k|x_{k-1})) z_{k-1}(x_k) - z_{k-1}(x_{k-1}) \right) e_{x_{k-1}}$
\ENDFOR
\end{algorithmic}

\end{algorithm}

The corresponding ODE (in the sense of Theorem~\ref{thm:convergencestochasticalgorithm} is given by
\begin{equation} \label{eq:ODEzlearning} \dot z(t) = - D (I-H) z(t),\end{equation}
where $D$ is a diagonal matrix given by $d_{ii} = q_i$, where $(q_i)$ denotes the invariant probability distribution of the Markov chain given by $[q_{ij}]$.

It is immediate from the Perron-Frobenius theorem that the eigenvalues of the matrix $I-H$ are strictly contained in the closed right halfplane, with a one-dimensional eigenspace corresponding to the zero eigenvalue and all other eigenvalues having strictly positive real part. This still holds for a multiplication of $I-H$ by an arbitrary diagonal matrix with positive diagonal entries, but this is less immediate (see e.g. \cite[Exercise 2.5.2]{hornjohnsontopics} for the nonsingular case). Therefore the linear subspace spanned by $z^*$ is globally attracting, and the Z-learning algorithm converges to this subspace by Theorem~\ref{thm:convergencestochasticalgorithm} (iii).

\subsection{D-stability as a necessary condition for convergence of stochastic approximation algorithms}
\label{sec:diagonalstabilitynecessary}
In the previous example, the positive stability of $I-H$ carried over to a multiplication by a positive diagonal matrix, $D(I-H)$, irrespective of the kind of diagonal matrix $D$. This kind of stability (invariant under left- (or right-)multiplication by an arbitrary positive diagonal matrix) is called D-stability in the literature (see e.g. \cite[Section 2.5]{hornjohnsontopics}), and the major difficulty with establishing local stability of the KL-learning algorithm consists of showing D-stability for the corresponding linearized ODE.

\section{Theoretical analysis of KL-learning}
\label{sec:convergence}

The KL-learning algorithm (Algorithm~\ref{alg:kllearning}) works well in practice, but for a rigorous theoretical analysis of its behaviour we need to make a few modifications, as given by Algorithm~\ref{alg:kllearninganalysis}. 

The modifications of Algorithm~\ref{alg:kllearninganalysis} with respect to Algorithm~\ref{alg:kllearning} are:
\begin{itemize}
\item[(i)] The values of $z$, $\lambda$ and $\Delta$ are indexed by the time parameter $k$ to keep track of all values;
\item[(ii)] Instead of a single step size $\gamma > 0$ and a finite time horizon $M \in \ints$ we consider an infinite time horizon and a decreasing sequence of stepsizes $(\gamma_k)$;
\item[(iii)] At every iteration, if necessary, a projection is performed (in the computation of $\Delta_k$) to ensure that $\lambda_k \geq \lambda_{\min} := \min_{i,j \in \{1, \hdots, n\}} \exp(-\beta c(j|i)) / 2$.
\end{itemize}

The modification (i) is purely a notational matter. Modification (ii) is standard in the analysis of stochastic approximation algorithms. If we would keep the stepsize constant the theoretical analysis would be harder. The practical effect of keeping the stepsize fixed is that the values of $(z_k, \lambda_k)$ will oscillate around the theoretical solution $(z^*, \lambda^*)$ with a bandwith depending on $\gamma$. Modification (iii) has minimal practical effect; we have not seen cases in which the projection step was actually made. The theoretical solution $\lambda^*$ satisfies $\lambda^* \geq 2 \lambda_{\min}$ by theory on nonnegative matrices \cite[Corollary 8.1.19]{hornjohnsonmatrixanalysis}. The constant 2 is arbitrary, chosen to ensure that $\lambda^*$ lies well above $\lambda_{\min}$. By Lemma~\ref{lem:invariants}, $\lambda_k$ is bounded from above, so there is no need to prevent $\lambda_k$ from growing large.

In this section we will write $\reals_+^n := \left\{ x \in \reals^n : x_i > 0 \ \mbox{for} \ i = 1, \hdots, n\right\}$.

In the discussion below we will only refer to Algorithm~\ref{alg:kllearninganalysis}.

\begin{algorithm}[ht]
\caption{KL-learning, notation for analysis}
\label{alg:kllearninganalysis}

\begin{algorithmic}
\STATE $\lambda_0 \gets 2 \lambda_{\min} = \min_{i,j=1,\hdots,n} \exp(-\beta c(j|i))$
\STATE $z_0(i) \gets \frac 1 n \lambda_0$, for all $i=1,\hdots, n$, 
\STATE $x_0 \gets \mbox{any state in $S$}$
\FOR{$k = 1 \ \TO \ \infty$}
\STATE $x_k \gets \mbox{independent draw from $q(\cdot|x_{k-1})$}$
\STATE $\Delta_k \gets \max \left\{ \exp(-\beta c(x_k|x_{k-1})) z_{k-1}(x_k)/ \lambda_{k-1} - z_{k-1}(x_{k-1}), (\lambda_{\min} - \lambda_{k-1})/\gamma_k \right\}$
\STATE $z_k \gets z_{k-1}$
\STATE $z_k(x_{k-1}) \gets z_{k-1}(x_{k-1}) + \gamma_k \Delta_k$
\STATE $\lambda_k \gets \lambda_{k-1} + \gamma_k \Delta_k$
\ENDFOR
\end{algorithmic}

\end{algorithm}

The initial value for $\lambda_0$ is moreless arbitrary, but it is important that $||z_0||_1 = \lambda_0$ and $\lambda_0 \in K$ with $K$ given by Lemma~\ref{lem:invariants}. We impose the following conditions.

\subsection{Hypothesis}
\label{hyp:kllearning}
\begin{itemize}
\item[(i)] Let $(\gamma_k)_{k \in \ints}$ be a sequence of nonnegative real numbers such that
$\sum_{m=1}^{\infty} \gamma_k = \infty$, $\sum_{m=1}^{\infty} \gamma_k^2 < \infty.$
\item[(ii)] Let $q(\cdot|\cdot)$ be irreducible aperiodic Markov transition probabilities on the finite state space $S = \{ 1, \hdots, n\}$ with invariant probabilities $q_i$, $i \in S$;
\item[(iii)] Let $c \in \reals^{n \times n}$;
\item[(iv)] Let the matrix $H = [h_{ij}] \in \reals^{n \times n}$ be given by $h_{ij} = \exp(-\beta c_{ij}) q_{ij}$ and the diagonal matrix $D = [d_i] \in \reals^{n \times n}$ by $d_i = q_i$;
\item[(iv)] Define continuous time processes $(z^k(t))_{t \geq 0}$ and $(\lambda^k(t))_{t \geq 0}$ for $k = 0, 1, \hdots$ as in Hypothesis~\ref{hyp:generalconvergence} (vi).
\end{itemize}

\subsection{Theorem (Convergence of KL-learning)}
\label{thm:convergencekllearning}
Consider Algorithm~\ref{alg:kllearninganalysis} under the conditions of Hypothesis~\ref{hyp:kllearning}. Then
\begin{itemize}
\item[(a)] With full probability, for any sequence of processes $(z^k, \lambda^k)$ there exists a subsequence uniformly on bounded intervals to some continuous functions $(z, \lambda)$, $z : [0,\infty) \rightarrow \reals^n_+$ and $\lambda : [0,\infty) \rightarrow (0,\infty)$.
\item[(b)] The trajectories of Algorithm~\ref{alg:kllearninganalysis}, as well as the limiting functions $(z,\lambda)$ given by (b), are constrained to a closed, bounded set $K$ given by~\eqref{eq:kllearningcompact}.
\item[(c)] Such a limit $(z, \lambda)$ satisfies the ODE
\begin{equation} \label{eq:ode_kllearning} \left\{ \begin{array}{ll}
\dot z(t) = f(z(t), \lambda(t)) + w, \\
\dot \lambda(t) = h(z(t), \lambda(t)) + \mu, \quad & t \geq 0, \end{array} \right.  \end{equation}
with $f : \reals^n_+ \times (0,\infty)\rightarrow \reals^n$ and $h : \reals^n_+ \times (0,\infty) \rightarrow \reals$ given by
\begin{equation}
\label{eq:odefunctions}
f(z, \lambda) := D \left(\frac 1 {\lambda} H - I \right) z, \quad h(z, \lambda) := \1^T f(z, \lambda),
\end{equation}
where
\[ D = \diag(q(1), \hdots, q(n)),\]
with $q$ the unique invariant probability distribution for the Markov chain with transition probabilities $q_{ij}$. Here $w \in \reals^n_{\geq 0}$ and $\mu \geq 0$ denote the minimum force necessary to keep $(z(t),\lambda(t))$ in $K$ (the continuous time equivalent of the projection step to ensure that $\lambda_k \geq \lambda_{\min}$; see \cite[Section 4.3]{kushneryin}).
\item[(d)] The ODE~\eqref{eq:ode_kllearning} admits a unique equilibrium $(z^*, \lambda^*)$ in the interior of $K$, where $H z^* = \lambda^* z^*$ and $||z^*||_1 = \lambda^*$.
\item[(e)] If any of the conditions of Proposition~\ref{prop:stabilitysufficientconditions} hold, then the equilibrium $(z^*,\lambda^*)$ as mentioned under (d) is locally asymptotically stable.
\end{itemize}

\emph{Proof.} By Lemma~\ref{lem:invariants} the trajectories of Algorithm~\ref{alg:kllearninganalysis} are constrained to the compact set $K$ given by~\eqref{eq:kllearningcompact}. We may apply a variant Theorem~\ref{thm:convergencestochasticalgorithm} suitable for projected algorithms (see \cite[Theorem 6.6.1]{kushneryin} to deduce (a), (b) and (c), where for (c) we may use Lemma~\ref{lem:odefunctions}. Results (d) and (e) follow from Propositions~\ref{prop:equilibrium}, Remark~\ref{remark:lowerdimension}, \ref{prop:equilibriumstability} and~\ref{prop:stabilitysufficientconditions}, where we note that no projection force is necessary in the interior of $K$.\qed

\subsection{Lemma (algorithm invariants)} 
\label{lem:invariants}
Under Hypothesis~\ref{hyp:kllearning}, the trajectories $(z_k,\lambda_k)$ of Algorithm~\ref{alg:kllearninganalysis} are contained in the compact set
\begin{equation} \label{eq:kllearningcompact} K = \left\{ (z, \lambda) \in [0,M]^n \times [\lambda_{\min}, n M] : ||z||_1 = \lambda \right\},\end{equation}
with $M := \max_{i,j=1, \hdots, n} \exp(-\beta c_{ij})$;

\emph{Proof.} Note that, by the maximum operation in the algorithm, $\Delta_k \geq \exp(-\beta c(x_k|x_{k-1})) z_{k-1}(x_k) - z_{k-1}(x_{k-1}))$. The update for $z_k$ therefore satisfies
\begin{equation} \label{eq:updatez} z_k(x_{k-1}) \geq (1-\gamma_k) z_{k-1}(x_{k-1}) + \gamma_k \exp(-\beta c(x_k|x_{k-1})) z_{k-1}(x_k)/\lambda_{k-1} > 0 \end{equation}
It follows immediately by induction that $z_k(i) > 0$ for $i =1, \hdots, n$ and $k = 1,\hdots, n$. The update-rule $\lambda_k \gets \lambda_{k-1} + \gamma_k \Delta_k$ for $\lambda_k$ ensures that $\lambda_k = ||z_k||_1 \geq \lambda_{\min}$ for all $k = 1, \hdots, n$.
We will show by induction that $z_k(i) \leq M$ for all $k \in \ints$ and $i =1, \hdots, n$. Recall that $z_0(i) \leq \frac M n \leq M$ for $i=1,\hdots, n$. Suppose $z_{k-1}(i) \leq M$ for all $i$, and some $k \in \ints$. If no projection occurs
\[ z_k(x_{k-1}) \leq (1-\gamma_k) z_{k-1}(x_{k-1}) + \gamma_k \max_{i,j=1,\hdots,n} \exp(-\beta c_{ij}) \leq M,\]
where we used that $z_{k-1}(i) \leq \lambda_{k-1}$ for all $i$. If projection does occur,
\[ z_k(x_{k-1}) \leq \lambda_k = \lambda_{k-1} + \gamma_k (\lambda_{\min} - \lambda_{k-1}) / \gamma_k = \lambda_{\min} < M.\]
\qed

\subsection{Lemma}
\label{lem:odefunctions}
The function $\overline g$ corresponding to Algorithm~\ref{alg:kllearninganalysis} in the sense of Theorem~\ref{thm:convergencestochasticalgorithm} is given by
\begin{align*}
\overline g(z, \lambda) = \begin{bmatrix} f(z,\lambda) \\ h(z, \lambda) \end{bmatrix},
\end{align*}
where $f$ and $h$ are given by~\eqref{eq:odefunctions}.

\emph{Proof.}
A straightforward computation.
\qed

\subsection{Proposition}
\label{prop:equilibrium}
Suppose $D$ is a diagonal matrix with positive entries and $H$ a nonnegative matrix.
Suppose $f : \reals^n \times \reals \rightarrow \reals^n$ and $h : \reals^n \times \reals \rightarrow \reals$ are given by~\eqref{eq:odefunctions}.

Consider the ODE~\eqref{eq:ode_kllearning}
with initial values $(z(0),\lambda)$ such that $z(0) > 0$ and $\lambda(0) = ||z(0)||_1$.

The orbits $(z(t), \lambda(t))$ satisfy $z(t) \geq 0$ and $\lambda(t) = ||z(t)||_1 > 0$ for all $t \geq 0$.
Furthermore a point $(z, \lambda) \in \reals^n \times \reals$ is an equilibrium point 
if and only if $H z = \lambda z$.

\emph{Proof.}
Suppose that for some $t \geq 0$, we have $z(t) \geq 0$, $\lambda(t) > 0$.
If for some component $z(t)(k)$ of $z(t)$ we have $z(t)(k) = 0$, then
\[ [f(z(t), \lambda(t))](k) = \left[ D \left(\frac 1 {\lambda(t)} H - I \right) z(t)\right](k) = \frac 1 { \lambda(t)} [D H z(t)](k) \geq 0.\]
Also, as $\lambda \downarrow 0$, then 
\[ h(z, \lambda) = \1^T D \left( \frac 1 {\lambda} H- I \right) z \rightarrow \infty,\]
so that $\lambda$ always remains positive.
For $z(t) \geq 0$, 
\[ \frac{d}{dt} ||z(t)||_1 = \sum_{k=1}^n \dot z(t)(k) = \1^T \dot z(t) = g(z(t),\lambda(t)) = \dot \lambda(t).\]

Because $z(0)(k) > 0$ for all $k = 1, \hdots, n$ and $\lambda(0) = ||z(0)||_1$, it follows that $z(t) \geq 0$  for all $t$ and $||z(t)|| = \lambda(t)$.

It is straightforward that $Hz = \lambda z$ if and only if $(z, \lambda)$ is an equilibrium point.
\qed

\subsection{Remark}
\label{remark:lowerdimension}
In light of the proposition above, we may consider the dynamical system 
\begin{equation}
\label{eq:ode2}
\left\{ \begin{array}{l} \dot z(t) = D \left( \frac 1 {||z(t)||_1} H - I \right) z(t), \\
         z(0) = z_0, \end{array} \right.
\end{equation}
with $z_0(k) > 0$ for all $k = 1, \hdots, n$, instead of~\eqref{eq:ode_kllearning}, thus reducing the dimensionality from $\reals^{n+1}$ to $\reals^n$.

\subsection{Definition}
A matrix $A \in \reals^{n \times n}$ is called \emph{stable} (or \emph{strictly stable}) if all eigenvalues $\lambda \in \sigma(A)$ satisfy $\Re \lambda \leq 0$ (or $\Re \lambda < 0$, respectively).

\subsection{Proposition}
\label{prop:equilibriumstability}
Suppose $H$ is a nonnegative irreducible matrix and $D$ is a diagonal matrix with positive entries on the diagonal. 

Then~\eqref{eq:ode2} has a unique equilibrium point $z^*$. This equilibrium point satisfies 
\begin{itemize}
\item[(i)] $||z^*||_1 = \lambda^* := \rho(H)$,
\item[(ii)] $z^* > 0$, and
\item[(iii)] $H z^* = \lambda^* z^*$.
\end{itemize}

The equilibrium is locally (asymptotically) stable if and only if the matrix $D(H- \lambda^* I - z^* \1^T)$ is (strictly) stable.

\emph{Proof.}
By Remark~\ref{remark:lowerdimension} above, we may apply Proposition~\ref{prop:equilibrium} to conclude that $z^*$ satisfies (iii) for some $\lambda^* > 0$, and $z^* \geq 0$, $z^* \neq 0$. Since $H$ is nonnegative and irreducible, there is, up to scaling by a positive constant, only a single eigenvector with nonnegative components. This is the Perron vector whose eigenvalue satisfies $\lambda^* =\rho(H)$, and which has only positive components, so that (i) and (ii) follow.

The linearization of $z \mapsto D \left( \frac 1 {||z||_1} H - I \right) z$ around $z$ is given by
\[ v \mapsto D \left( \frac 1 {||z||_1} H - I \right) v - \frac 1  {||z||_1^2} DH z \1^T v,\]
which reduces to
\[ v \mapsto D \left( \frac 1 {\lambda^*} (H - z^* \1^T) - I \right) v \] for $z = z^*$. Multiplication by $\lambda^*$ does not affect the stability properties of this matrix, so the stability of the equilibrium $z^*$ is determined by the spectrum of the matrix $D(H - z^* \1^T - \lambda^* I)$.
\qed

The stability of the matrix $D(H - \lambda^* I - z^* \1^T )$ seems to be a non-trivial issue. In Proposition~\ref{prop:stabilitysufficientconditions} below we collect some facts that we have already obtained. For this we need a lemma.

\subsection{Lemma}
\label{lem:nonsingular}
Suppose $H$ and $D$ satisfy the conditions of Proposition~\ref{prop:equilibriumstability}. Then the matrix $D(H-\lambda^* I - z^* \1^T)$, with $\lambda^* = \rho(H)$ and $z^*$ the corresponding positive eigenvector, is nonsingular.

\emph{Proof.} 
We will omit *-superscripts in this proof, so $z = z^*$ and $\lambda = \lambda^*$. Write $A = H - \lambda I - z \1^T$.

Let $w$ and $z$ denote the left and right Perron-Frobenius eigenvectors of $H$. Note that $(z,w) > 0$ and $(\1, z) > 0$. Also note that any $\zeta \in \reals^n$ may be written as $\zeta = \alpha z + \eta$, with $\eta \perp \1$, by picking $\alpha = (\zeta, \1) / (z, \1)$ and $\eta = \zeta - \alpha z$. Therefore we may choose a basis of $\reals^n$ consisting of the vector $v_1 = z$ and some vectors $v_2, \hdots, v_n$ spanning $\1^{\perp}$. Let $S$ denote the matrix with columns $v_1, \hdots, v_n$. Then the first column of $S^{-1} (H - \lambda I) S$ consists of zeroes since $(H - \lambda I) z = 0$, and only the first column of $S^{-1} z \1^T S$ is nonzero. So adding $S^{-1} z \1^T S$ to $S^{-1}(H -\lambda I) S$ only increases the range of the resulting matrix. Therefore, $\mathrm{range} (H - \lambda I) \subset \mathrm{range}(A)$.

Since $w^T (H - \lambda I) = 0^T$, we have that $w$ is perpendicular to the range of $H - \lambda I$. But $w$ is not perpendicular to the range of $A$ since $(w,z) > 0$. In other words, the inclusion $\mathrm{range}(H - \lambda I) \subset \mathrm{range} (A)$ is strict, $\mathrm{rank}(A) > \mathrm{rank}(H - \lambda I) = n-1$, so that $\mathrm{rank}(A) = n$ and $\det(A) \neq 0$. Hence also $\det ( DA ) \neq 0$.
\qed

\subsection{Proposition}
\label{prop:stabilitysufficientconditions}
Suppose $H$ and $D$ satisfy the conditions of Proposition~\ref{prop:equilibriumstability}. The matrix $D(H - \lambda^* I - z^* \1^T )$, with $\lambda^* = \rho(H)$ and $z^*$ the corresponding positive eigenvector, is strictly stable in any of the following cases:
\begin{itemize}
\item[(i)] $D = \beta I$ for some $\beta > 0$,
\item[(ii)] $\1^T H = \lambda^* \1^T$ (so $\1^T$ is a left Perron vector),
\item[(ii)] $D, H \in \reals^{2 \times 2}$.
\end{itemize}

\emph{Proof.}
As before, we will omit *-superscripts in this proof, so $z = z^*$ and $\lambda = \lambda^*$.

\begin{itemize}
\item[(i)] 

Suppose $v$ is an eigenvector of $A = H - \lambda I - z \1^T$ with eigenvalue $\mu \neq 0$. Define $w = v + \left( \frac{ \1^T v}{\mu} \right) z$. Then
\[ (H - \lambda I) w = (H - \lambda I) v = \left( \mu v + z \1^T v  \right) = \mu w,
 \]
which shows that $\mu \in \sigma(H - \lambda I)$, and since $\mu \neq 0$, it follows that $\mathrm{Re} \ \mu < 0$. So all $\mu \in \sigma(A)$ have $\mathrm{Re} \ \mu < 0$, except possibly the case where $\mu = 0$ but this case is excluded by Lemma~\ref{lem:nonsingular} above.
\item[(ii)]

Let $B = (H - \lambda I) \diag(z)$. Then $B$ has a positive diagonal and negative off-diagonal entries. Furthermore $B \1 = 0$ and $\1^T B = 0$. This shows that $B$ is row diagonally dominant and column diagonally dominant, so that the same holds for $B + B^T$. It follows that $B + B^T$ is positive semidefinite. Note that $B + B^T$ is a singular $M$-matrix (see \cite{hornjohnsontopics}, Section 2.5.5). Since $H$ is irreducible, also $B + B^T$ is irreducible. Therefore the nullspace of $B + B^T$ is one-dimensional and it is spanned by $\1$. Also $z \1^T \diag(z) = z z^T$ is symmetric positive semidefinite, and $(\1, z z^T \1) > 0$. It follows that $B + B^T + 2 z z^T$ is positive definite. Multiplying on both sides by $D$ (a congruence transform) gives that $D (H-\lambda I - z \1^T) \diag(z) D + \diag(z) D (H -\lambda I -z \1^T)^T D$ is symmetric positive definite. By Lyapunov's theorem \cite[Theorem 2.2.1]{hornjohnsontopics}, it follows that $D(H - \lambda I - z \1^T)$ is strictly stable.

\item[(iii)]
By (i) we have that $A = H - \lambda I -z \1^T$ is strictly stable. In 2 dimensions, this is equivalent to $\det(A) > 0$ and $\mbox{tr}(A) <0$. This immediately implies that $\det (DA) = \det(D) \det(A) > 0$. We compute
\[ \mbox{tr}(D A) = \mbox{tr} \begin{bmatrix} d_1 & 0 \\ 0 & d_2 \end{bmatrix} \begin{bmatrix} h_{11} - \lambda - z_1 & h_{12} - z_1 \\ h_{21} - z_2 & h_{22} - \lambda - z_2 \end{bmatrix} = d_1 (h_{11} - \lambda - z_1) + d_2 (h_{22} - \lambda - z_2) < 0,\]
since the diagonal of $H - \lambda I$ has nonpositive entries (which may be seen by \cite{hornjohnsonmatrixanalysis}, Theorem 8.3.2).
\end{itemize} \qed

We think the above proposition can be generalized significantly. In fact, we propose the following conjecture. If the conjecture holds, Theorem~\ref{thm:convergencekllearning} (e) can be formulated unconditionally.

\subsection{Conjecture}
\label{conj:stability}
Suppose $H$ and $D$ satisfy the conditions of Proposition~\ref{prop:equilibriumstability}. Then the matrix $D(H -\lambda^*I - z^* \1^T )$ is strictly stable.

% We are preparing a theoretical analysis of more general dynamical systems that have the Perron vector as equilibrium, including a more detailed analysis of their stability properties, see \cite{bierkenskappen_mathematical_2012}.

\section{Numerical experiment}

\label{sec:numerical}
\begin{figure}
\begin{center}
\subfigure[Grid world]{\includegraphics[width=0.45 \linewidth]{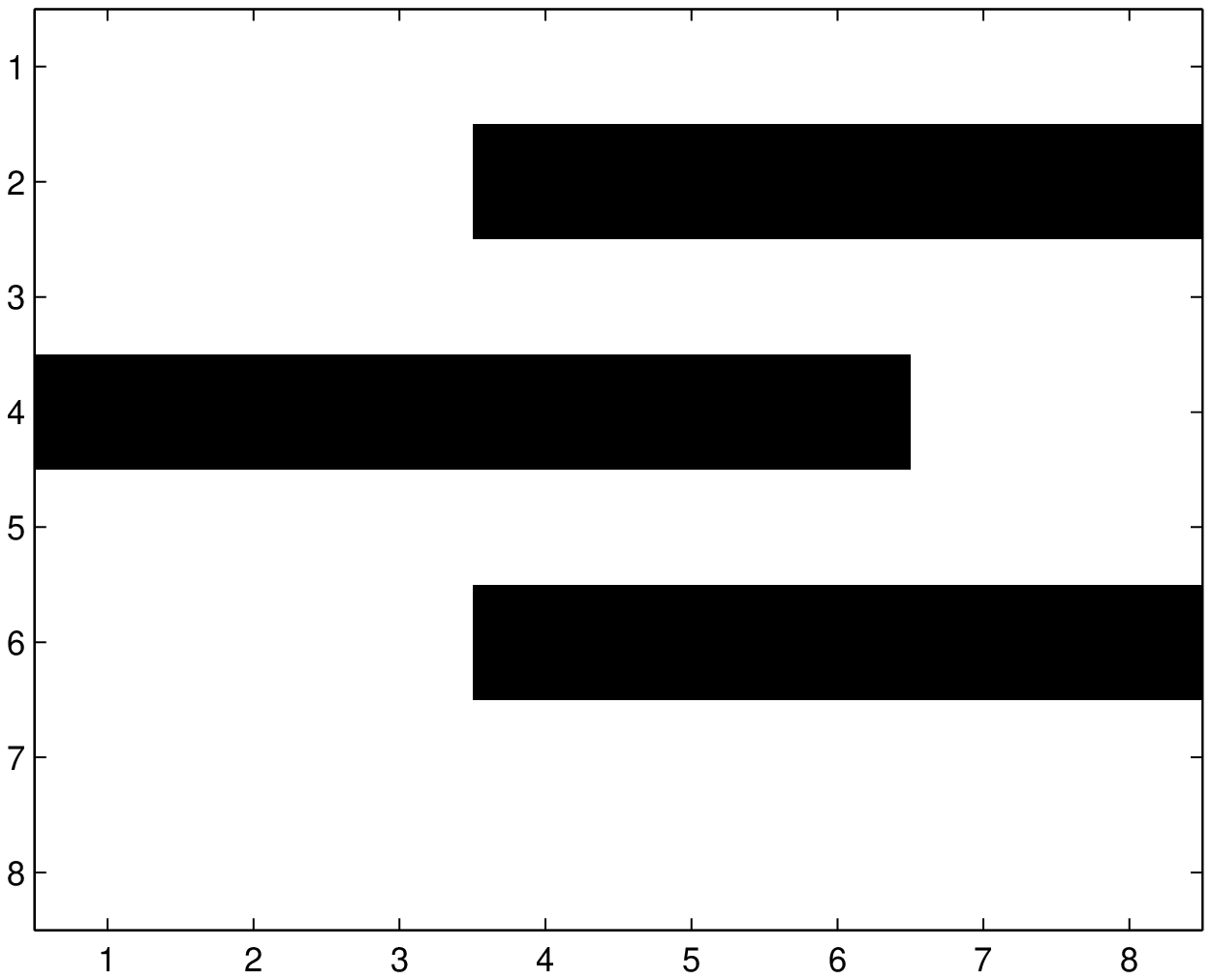}}  
\subfigure[Value function]{\includegraphics[width=0.45 \linewidth]{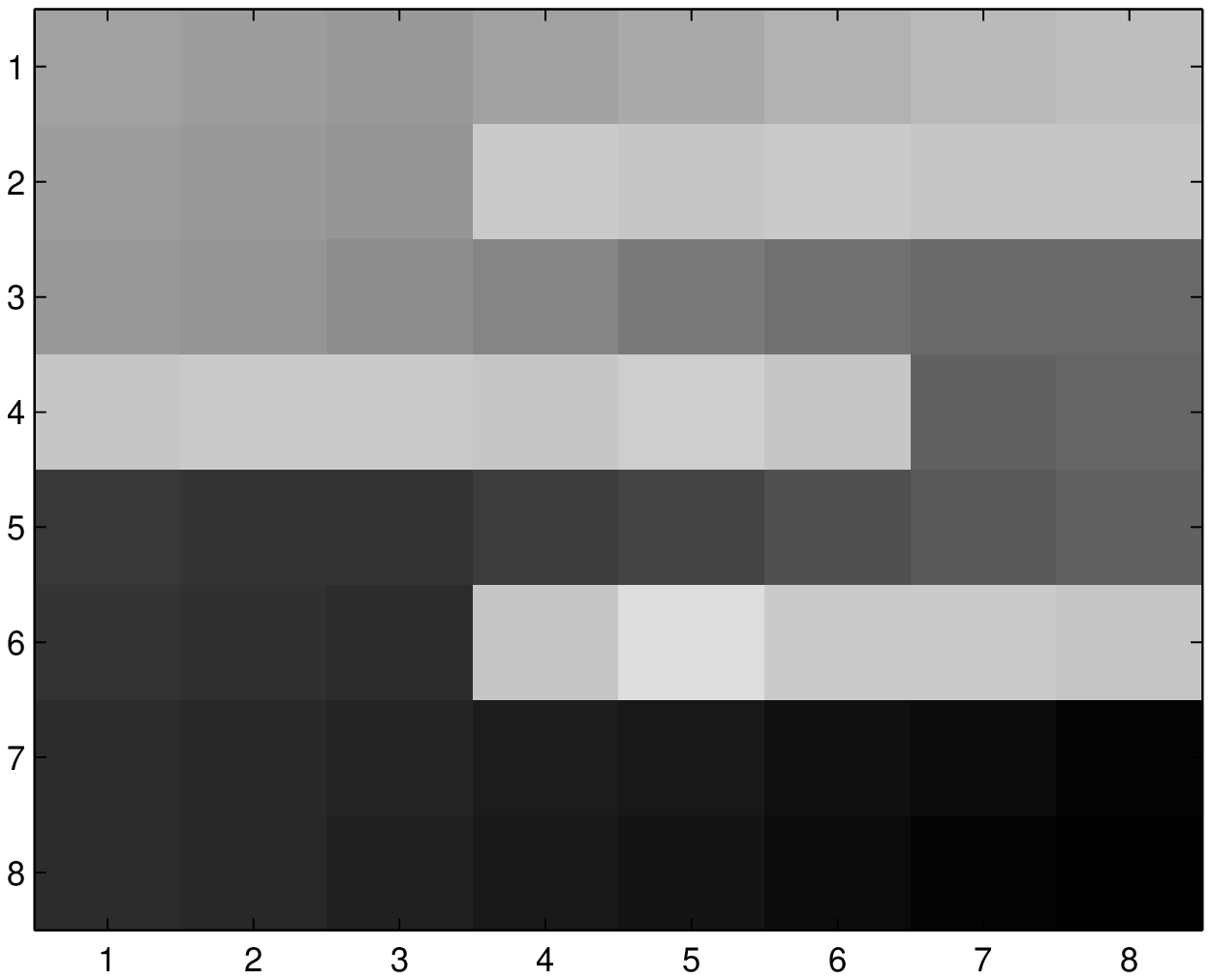}} \\
\subfigure[Comparison between Z-learning and KL-learning]{\includegraphics[width=0.45 \linewidth]{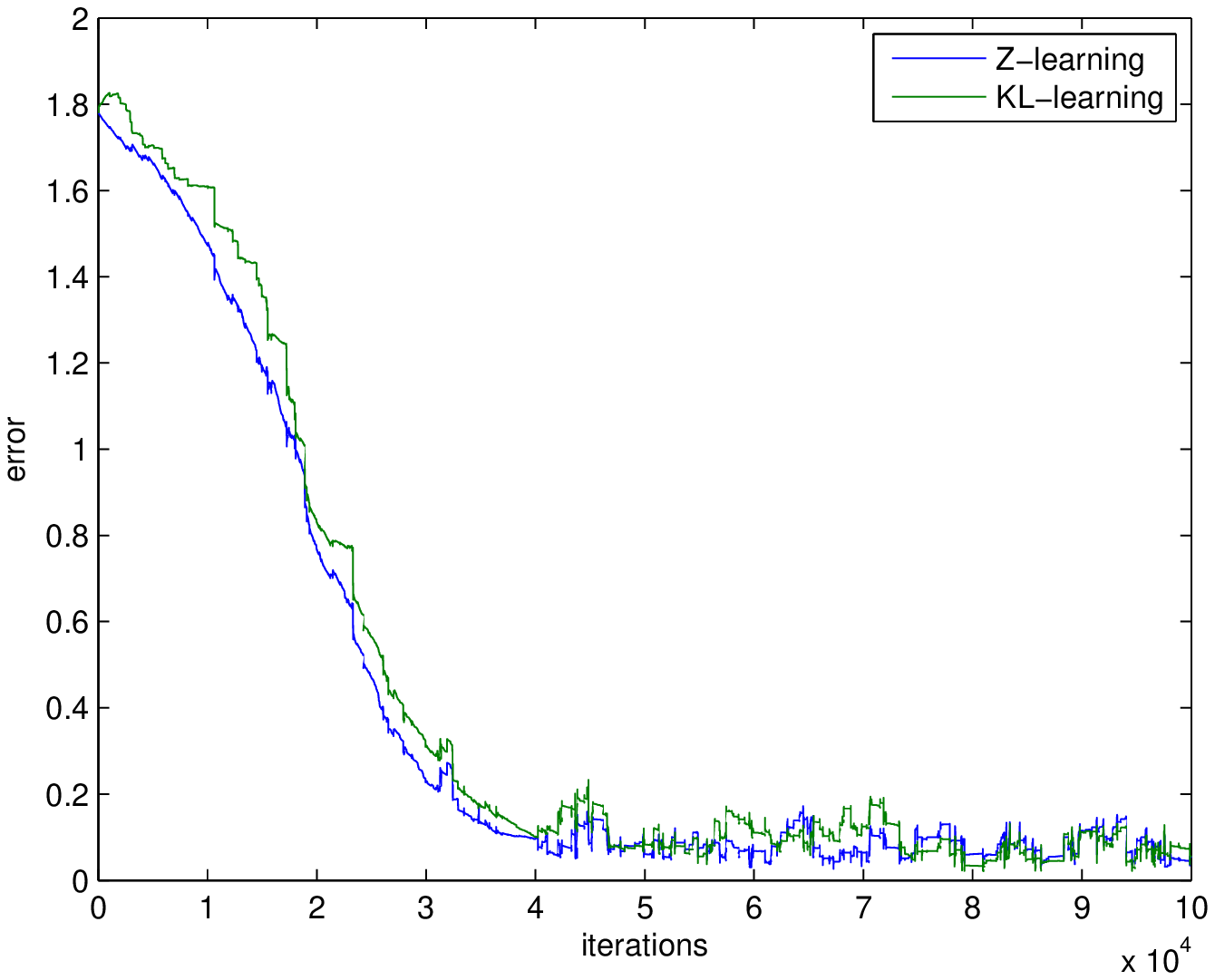}}
\subfigure[Power method]{\includegraphics[width=0.5 \linewidth]{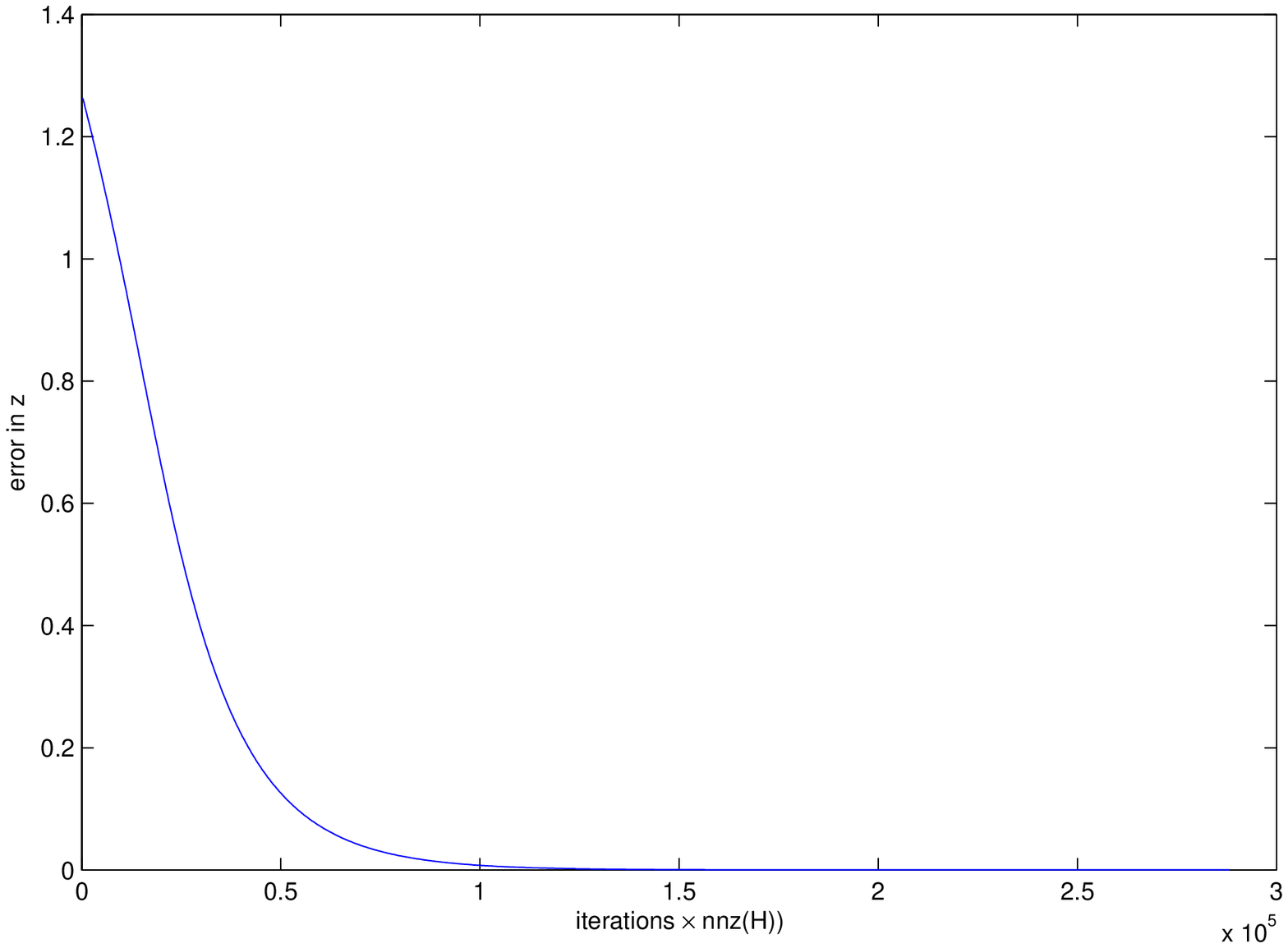}}
\caption{Numerical experiment}
\label{fig:experiment}
\end{center}
\end{figure}

Consider the example of a gridworld (Figure~\ref{fig:experiment} (a)), where some walls are present in a finite grid. Suppose the uncontrolled dynamics $q$ allow to move through the walls, but walking through a wall is very costly, say a cost of 100 per step through a wall is incurred. Where there is no wall, a cost of 1 per step is incurred. There is a single state, in the bottom right, where no costs are incurred. The uncontrolled dynamics are such that with equal probability we may move left, right, up, down or stay where we are (but it is impossible to move out of the gridworld). The value function for this problem can be seen in Figure~\ref{fig:experiment} (b). In order to be able to compare our algorithm to the original Z-learning algorithm, the cost vector is normalized in such a way that $\lambda^* = 1$, so that Z-learning converges on the given input.

The result of running the stochastic approximation algorithm, with a constant gain of $\gamma = 0.05$ is portrayed in Figure~\ref{fig:experiment} (c), where it is compared to Z-learning (see Section~\ref{sec:zlearning} and \cite{todorov2007}). This result may also be compared to the use of the power method in Figure~\ref{fig:experiment} (d). Here the following version of the power method is used, in order to be able to give a fair comparison with our stochastic method.
\[ z_k = z_{k-1} + \gamma_k(H z_{k-1} - z_{k-1}).\] Note that for each iteration, the number of operations is (for sparse $H$) proportional to the number of non-zero elements in $H$. In the stochastic method the number of operations per iteration is of order 1. Comparing the graphs in Figure~\ref{fig:experiment} (c) and (d), we see that KL-learning does not disappoint in terms of speed of convergence, with respect to Z-learning as well as the power method.

\section{Discussion}
\label{sec:discussion}

The strength of KL control is its very general applicability. The only requirements are the existence of some uncontrolled dynamics governed by a Markov chain, and some state or transition dependent cost. The Markov chain may actually be derived from a graph of allowed transitions, giving every allowed transition equal probability. A disadvantage is that we cannot directly influence the control cost; it is determined by the KL divergence.

KL control is very useful if we know which moves (e.g. in a game) are allowed and we wish to find out which moves are best. The control cost of KL divergence form has a regularizing effect: no move will be made with probability one (unless it is the only allowed move). You could say that there is always a possibility to perform an exploratory move, instead of an exploiting move, under the controlled dynamics.

This immediately suggests the use of KL learning as a reinforcement learning algorithm. The initial transition probabilities represent exploratory dynamics. At every iteration, we could compute a new version of the optimal transition probabilities and use these as a new mixture of exploitation and exploration. The practical implications of this idea will be the topic of further research.

The KL learning algorithm seems to work well in practice and a basis has been provided for its theoretical analysis. Some questions remain to be answered. In particular, if Conjecture~\ref{conj:stability} is true, then regardless of the structure of the problem we know that the solution of the control problem is a locally asymptotically stable equilibrium of the algorithm. It would be even more convenient if a Lyapunov function for the ODE~\eqref{eq:ode_kllearning} could be found, which would imply global convergence of KL learning.

So far numerical results indicate that KL learning is a reliable algorithm. In the near future we will apply it to practical examples and evaluate its performance relative to other reinforcement learning algorithms.

\end{document}